\documentclass{amsart}
\usepackage{amssymb,amsfonts,latexsym}

\renewcommand{\Bbb}{\mathbb}
\renewcommand{\frak}{\mathfrak}

\theoremstyle{plain}
\newtheorem{thm}{Theorem}[section]
\newtheorem{prop}[thm]{Proposition}
\newtheorem{lem}[thm]{Lemma}
\newtheorem{cor}[thm]{Corollary}
\newtheorem{lemma}[thm]{Lemma}
\theoremstyle{definition}
\newtheorem{remark}[thm]{Remark}

\title[Liv\v{s}ic's Theorem, Superrigidity, and Anosov Actions]
{On Liv\v{s}ic's Theorem, Superrigidity, and Anosov Actions of 
Semisimple Lie Groups}
\author{Edward R. Goetze}
\thanks{The first author was supported in part by grants from the NSF and the 
University of Michigan.}
\address{University of Michigan,
Department of Mathematics,
Ann Arbor, MI 48109-1092}
\email{goetze@math.lsa.umich.edu}

\author{Ralf J. Spatzier}
\thanks{The second author was supported in part by a grant from the NSF}
\address{University of Michigan,
Department of Mathematics,
Ann Arbor, MI 48109-1092}
\email{spatzier@math.lsa.umich.edu}

\begin{document}

\begin{abstract}
We prove a generalization of Liv\v{s}ic's Theorem on the vanishing
of the cohomology of certain types of dynamical systems.
As a consequence, we strengthen a result due to Zimmer concerning
algebraic hulls of Anosov actions of semisimple Lie groups.  Combining
this with Topological Superrigidity, we find a H\"{o}lder geometric
structure for multiplicity free Anosov actions.
\end{abstract}

\bibliographystyle{plain}
\maketitle

\begin{center}
January 16, 1996
\end{center}

\section{Introduction}

During the last decade, Anosov actions of semisimple Lie groups
and their lattices have become a focal point in the study of rigidity 
properties of such groups.  Most importantly, local smooth rigidity 
has been established 
for various standard algebraic Anosov actions by a number of authors 
\cite{lewis1,hurder,kl1,kl2,klz,qian1,qian2,qian3}. 
The conjecture arose that all such actions are essentially
$C^{\infty}$-conjugate to algebraic actions. The proof of this conjecture 
for the special class of Cartan actions is the goal of
the current paper and its sequel \cite{gs2}.

In this paper, we introduce an additional geometric structure for certain
types of 
volume preserving Anosov actions of
a connected higher rank semisimple Lie group $G$ of the noncompact type on
a closed manifold $M$.  More precisely, 
we will find a H\"{o}lder Riemannian metric, a H\"{o}lder splitting $\bigoplus
E_{\lambda}$ of the
tangent bundle and a finite dimensional representation $\pi$ of $G$
such that the elements in a Cartan subgroup of $G$ expand and contract 
vectors in $E_{\lambda}$ {\em precisely} according to a weight of $\pi$. 
Note that Zimmer's superrigidity theorem for cocycles yields the same
conclusion with respect to a measurable Riemannian metric. The difference in
regularity however is crucial to our classification of such actions.

One fundamental tool required to obtain these geometric results is
a generalization of celebrated work of Liv\v{s}ic on the
cohomology of Anosov systems. 
Liv\v{s}ic showed in particular that an 
$\Bbb R$-valued H\"{o}lder cocycle which 
is measurably cohomologous to the trivial cocycle is 
H\"{o}lder cohomologous to the trivial cocycle \cite{livsic1}.
Liv\v{s}ic also obtained results for cocycles taking 
values in nonabelian Lie groups under the additional assumption 
that the cocycle evaluated on generators takes 
values sufficiently close to the identity.
Theorem \ref{livthm} and its 
corollaries 
provide bundle theoretic 
versions of Liv\v{s}ic's theorem as we will now explain. 
Suppose a group $G$ acts on a principal $H$-bundle $P \rightarrow M$ via bundle
automorphisms. Since any bundle has a measurable section, such an action is
measurably isomorphic to a skew product action on $M \times H$ via a cocycle 
$\alpha:G \times M \to H$.  Different trivializations of $P$ 
correspond to distinct yet
cohomologous cocycles. Note that $\alpha$ is 
measurably cohomologous to the trivial cocycle precisely
when there is a measurable $G$-invariant section of $P$. More generally,
there is a complete correspondence in the measurable category between bundle
theoretic and cohomological statements. This correspondence breaks down in the
continuous or smooth category. Indeed, $P$ may not admit any continuous
sections, and thus the $G$ action on $P$ may not give rise to a cocycle.
Nevertheless, Liv\v{s}ic's theorem generalizes. In particular, if there is a
measurable invariant section of $P$ and $H$ is the semidirect product
of compact and abelian groups, then there is also an invariant H\"{o}lder
section. 
%objective in this paper, it also has
%independent interest. In particular, it yields new applications to skew 
%products. We pursue this elsewhere \cite{gs3}.

Our primary application of this bundle theoretic 
version of Liv\v{s}ic's theorem is a description of the H\"{o}lder
algebraic hull of certain $G$ actions on the frame bundle,
where $G$ is a semisimple Lie group of higher rank without
compact factors. Algebraic hulls are invariants of the action, defined
in the measurable, continuous, H\"{o}lder and smooth categories for any
principal bundle action as above as long as the structure group $H$ is an
algebraic group. They are the smallest algebraic subgroup of $H$ for which
there is a reduction of $P$ of the relevant regularity which is invariant 
under the $G$-action. In the measurable case, a trivial algebraic hull is
equivalent to the associated cocycle being cohomologically trivial. 
In \cite{z10}, Zimmer was able to show that
the measurable algebraic hull of a higher rank semisimple Lie group
is reductive with compact center.  In Theorem \ref{THM2}, we strengthen
this description to include the H\"{o}lder algebraic hull of 
the frame bundle for Anosov actions
of higher rank semisimple Lie groups.
The proof uses our version of Liv\v{s}ic's theorem as well as arguments
from finite dimensional representation theory.
To produce the special H\"{o}lder Riemannian metric for 
multiplicity free Anosov actions of a semisimple
group, we combine our description of algebraic hulls with Zimmer's Topological 
Superrigidity Theorem \cite{erg1,z9}. The latter theorem allows one to
find continuous, H\"{o}lder or smooth sections of $P$ which transform under the
$G$ action essentially according to a finite dimensional representation $\pi$
of $G$ provided one can find sections of the same regularity 
of certain associated bundles to $P$
invariant under a parabolic subgroup of $G$, e.g. Grassmann bundles. 
In our case, we can use the stable distribution of an 
Anosov element of $G$ as the 
H\"{o}lder section of a Grassmann bundle. 

We wish to thank Gopal Prasad 
for several enlightening conversations and Nantian Qian for pointing out
a gap in a preliminary version of this paper.

\section{Liv\v{s}ic's Theorem and the H\"{o}lder Algebraic Hull}

The main purpose of this section is to prove our generalization of
Liv\v{s}ic's theorem, Theorem \ref{livthm}.  Let us we begin
by mentioning some of the basic notions central to our presentation.

\subsection{Preliminaries}
Let $G$ be a Lie group, possibly discrete. Suppose $G$ 
acts smoothly (at least $C^{1+\theta})$ and locally freely on a manifold 
$M$ with a norm $\parallel \cdot \parallel$ associated to a
Riemannian metric. Call an element 
$g \in G$ {\em regular}  or {\em normally hyperbolic} 
if there exist a continuous decomposition of the tangent bundle
\[TM = \tilde{E}^+ _g + \tilde{E}^0  + \tilde{E}^- _g\]
into $g$-invariant subbundles,  and positive constants 
$\tilde{C}>1$, $\tilde{A},\tilde{B}> 0$ such
that for every $m \in M$, for every positive integer $n$, and for
every $v \in \tilde{E}_g^\pm(m)$,
\begin{equation*}
\frac{1}{\tilde{C}} \|v\|e^{-n\tilde{B}} \le \|Tg^{\mp n}v\| \le \tilde{C}
\|v\|e^{-n\tilde{A}},
\end{equation*}
and such that $E^0$ is the tangent distribution of the $G$-orbits.
Call a $G$ action {\em Anosov}  or {\em normally hyperbolic} 
if it contains a normally  hyperbolic element.
While $\tilde{E}^-_g$ is clearly contained in the stable 
distribution $E^-_g$ of $g$,
$E^-_g$ may also contain elements tangent to the orbits of $G$. To analyze
this contribution, let  $q$ be the dimension of $G$. Consider the 
$q$-frames $\mathcal{Q}$ of
$M$ tangent to the $G$-orbits. Fix a basis $X_1, \ldots, X _q$ of $\frak{g}$.
Since the $G$-action is locally free, evaluating the $X_i$ defines a 
smooth section $\theta$ of $\mathcal{Q}$. In fact, $\mathcal{Q}$ is a 
trivial bundle with fiber the full frames of 
$\frak{g}$, on which $G$ acts via the adjoint representation. 
Note that $\theta$ transforms according to the adjoint transformation for $h
\in G$,
\[ Th (\theta(m)) = Ad(h) \theta (hm),\]
where $Th$ denotes the derivative of $h$.

For $g$ a normally hyperbolic element of the $G$ action as above, 
let $\mathcal{O}_g^-(m)$ correspond to the sum of the generalized 
eigenspaces of $Ad(g)$ of eigenvalue of modulus less than 1. Then 
\[ E^-_g = \tilde{E}^-_g \oplus\mathcal{O}_g^- \]
is precisely the stable distribution of $g$ on $M$. Since $\tilde{E}^-_g$ is
continuous and $\mathcal{O}_g^-(m)$ is smooth, $E^-_g$ is continuous.
Similarly, the unstable distribution $E^+_g $ is continuous and
splits as a sum  
\[ E^+_g = \tilde{E}^+_g \oplus \mathcal{O}_g^+ .\]
The neutral distribution $E^0 _g$ corresponds to the sum of the
generalized eigenspaces of $Ad (g)$ with eigenvalue of modulus 1. Since the
latter form a Lie subalgebra of $\frak{g}$, it follows 
that $E^0_g$ is an integrable distribution. In particular, $g$ is 
normally hyperbolic to the orbit foliation of a subgroup of $G$.
Finally, there are constants $C>1$, $A,B>0$ such
that for every $m \in M$, for every positive integer $n$, and for
every $v \in E_g^\pm(m)$,
\begin{equation}
\frac{1}{C} \|v\|e^{-nB} \le \|Tg^{\mp n}v\| \le C\|v\|e^{-nA}.
\label{anosoveqn}
\end{equation}

If $M$ is compact, these notions do not depend 
on the ambient Riemannian metric. 
Note that the splitting and the  constants in the definition above 
depend on the normally hyperbolic element. 
It is well known that the distributions $E^+_g$ and $E^-_g$ are integrable
and are tangent to $W^s_g$ and $W^u_g$, the stable and unstable foliations
of $g$.  In particular, $W^-_g(x)$ denotes the stable manifold through
$x \in M$.  This is a H\"{o}lder foliation whose leaves
are smoothly immersed submanifolds of $M$ \cite{anosov1,hps,mane1}.
Suppose additionally that an Anosov action preserves a volume. 
Then the stable and unstable foliations are absolutely continuous. 
This follows from \cite[Theorem 2.1]{ps1} since the normally 
hyperbolic element $g \in G$ expands the neutral distribution 
$E^0 _g$ subexponentially by our analysis above. It follows from 
the usual Hopf argument and the fact that the neutral foliation 
is contained in the $G$-orbits that the $G$-action is ergodic. 

Let $P \to M$ be a principal bundle over $M$ with structure group
$H$, a Lie group, and suppose $G$ 
is a group.  We say that $G$ acts on $P$ via {\em bundle automorphisms}
if the $G$ action on $P$ factors to a $G$ action on $M$.
Fix an $H$-invariant metric on $P$. Then we say that $G$ acts via
{\em H\"{o}lder bundle automorphisms} if the $G$ action on $P$
projects to a smooth (at least $C^{1+\theta}$) 
action on $M$ and if each element of $G$ is a H\"{o}lder homeomorphism 
of $P$.

%Let us discuss the relationship with the standard version of
%Liv\v{s}ic's Theorem.  
Let us discuss the relationship between bundle automorphisms and
cocycles.
First, we recall the notion of a cocycle.
Suppose that $G$ is a group which acts on a manifold $M$, and $H$ is a
Lie group.  A function
$\alpha: G \times M \to H$ is called a {\em cocycle} if it satisfies the
cocycle identity:
\[ \alpha(g_1g_2,m) = \alpha(g_1,g_2m)\alpha(g_2,m) \text{ for all }
g_1, g_2 \in G, \text{ and } m \in M. \]
When dealing with measurable cocycles, we require that the cocycle identity
hold only for almost all $m \in M$.
Two cocycles $\alpha,\beta:G \times M \to H$ are called {\em equivalent}
or {\em cohomologous} if there exists a function $\phi:M \to H$ such that
$\phi(gm)^{-1} \alpha(g,m) \phi(m) = \beta(g,m)$.  The regularity of
$\phi$ determines the regularity of the equivalence.  In particular,
we say that a cocycle is measurably (H\"{o}lder, smoothly, etc.) {\em trivial}
if it is measurably (H\"{o}lder, smoothly, etc.) equivalent to the
trivial cocycle.

Given a cocycle $\alpha: G \times M \to H$, one can construct
a $G$ action on the trivial bundle $P = M \times H$
via bundle automorphisms by defining
\[ g(m,h) = (mg, \alpha(g,m)h). \]
Note that the regularity of the $G$ action on $P$ will be the same as
that of the cocycle $\alpha$.  In particular, if $\alpha$ is
H\"{o}lder, then the $G$ action on $P$ will be via
H\"{o}lder bundle automorphisms.

Conversely, if $G$ acts on $P$ via bundle automorphisms, then with respect
to any trivialization of $P$ (e.g., a measurable one), one can construct
a cocycle which describes this action. In particular, if $\sigma:M \to P$ 
is a section, then there exists a cocycle $\alpha:G \times M \to H$ so
that $g \sigma(m) = \sigma(gm)\alpha(g,m)$.  We refer to $\alpha$ as 
the {\em cocycle corresponding to $\sigma$}.
Although different trivializations
yield different yet cohomologous cocycles, we shall abuse
notation and refer to this class
of cocycles as the cocycle for the $G$ action on $P$.

We now state our main result, and then list a number of
consequences.

\begin{thm}
\label{livthm}
Let $P \to M$ be a principal $H$ bundle over a compact connected manifold
$M$ where $H=\mathcal{K}\ltimes\mathcal{A}$, the semidirect product of
a compact group with an abelian group.  Suppose $G$ is a Lie group that acts 
via H\"{o}lder bundle automorphisms on $P$ such that 
the $G$ action on $M$ is Anosov with $a\in G$ normally hyperbolic.
Let $V$ be a transitive left $H$ space which admits an 
$H$-invariant metric.  
Then any $G$-invariant measurable section of the 
associated bundle $E_V \to M$ is actually H\"{o}lder.
\end{thm}
\begin{remark}
Just as in the case of Li\v{v}sic's original work, there exist $C^1$ and
$C^\infty$ versions of this result.  We will pursue this elsewhere \cite{gs3}.
\end{remark}
\begin{remark}
A key element in the proof lies in estimating the effects of conjugation
by elements in $H$, cf. Equation \ref{FF1}.  We remark that this theorem 
will hold whenever an appropriate bound on conjugation can be made.  
As Liv\v{s}ic does in his work, if we assume that the cocycle corresponding
to some bounded section takes
values sufficiently close to the identity, then it is possible to 
produce such a bound.  This yields the following result.
\end{remark}
\begin{cor}
Let $P \to M$ be a principal $H$ bundle over a compact connected manifold
$M$ where $H$ is any algebraic group. Suppose $G$ is a Lie group that acts
via H\"{o}lder bundle automorphisms on $P$ such that
the $G$ action on $M$ is Anosov with $a\in G$ normally hyperbolic.
\begin{enumerate}
\item[($\star$)]
Assume there exists a measurable section $\sigma:M \to P$ taking values in
a compact subset $K \subset P$
so that for the corresponding cocycle $\alpha:G \times M\to H$, for some
choice of inner product on $\frak{h}$, and for some sufficiently small
$\eta >0$, the operator norm of $Ad(\alpha(a,x))$ is less than $1+\eta$.
\end{enumerate}
Let $V$ be a transitive left $H$ space which admits an
$H$-invariant metric.
Then any $G$-invariant measurable section of the
associated bundle $E_V \to M$ is actually H\"{o}lder.
\end{cor}

Call a set $U \subset M$ a {\em H\"{o}lder (measurable)
generic} set if $U$ contains an open dense (conull) set.
Now suppose $G$ acts by automorphisms of a principal $H$-bundle
$P \rightarrow M$ where $H$ is an algebraic group.  Assume that the
$G$ action on $M$ is ergodic with respect to some invariant measure
$\mu$.
An algebraic subgroup $L \subset H$ is called a {\em H\"{o}lder (measurable)
algebraic hull} for the $G$ action on $P$ if
\begin{enumerate}
\item[(a)] there exist a $G$-invariant H\"{o}lder (measurable)
generic set $U \subset M$ and a $G$-invariant H\"{o}lder (measurable)
section of $E_{H/L} \mid _U \rightarrow U$, and
\item[(b)] the first assertion is false for any proper algebraic
subgroup of $L$.
\end{enumerate}

Algebraic hulls exist and are unique up to conjugacy.
See \cite{erg1,z9,z10} for discussions of the
properties of algebraic hulls and some of their geometric consequences.
In a sense, the algebraic hull describes how much of $H$ is involved in the
$G$ action.  We remark that we can also define these terms in the $C^r$ and
Lipschitz categories in the obvious fashion.

Of course, if $H$ admits a bi-invariant metric, then the homogeneous space
$H/H_1$ for $H_1 \subset H$ admits a left invariant metric.  
Thus, Theorem \ref{livthm} has an immediate application for algebraic hulls.

\begin{cor}
Let $P \to M$ be a principal $H$ bundle over a compact connected manifold
$M$ where $H \subset GL(n,\Bbb R)$ is algebraic.
Suppose $G$ acts via bundle automorphism on $P$ such
that the $G$ action on $M$ is Anosov.
If $H$ admits a bi-invariant metric, then
the measurable algebraic hull for the $G$ action on
$P$ and the H\"{o}lder algebraic hull for the $G$ action on $P$
are equal.
\label{livcor}
\end{cor}

\subsection{Proof of Theorem \ref{livthm}}
The basic strategy of the proof is similar to that of Liv\v{s}ic.
A significant difference in our proof is that,
unlike in Liv\v{s}ic's situation, we do not have a globally defined section
of $P$ giving rise to a cocycle.  As it turns out, it is sufficient to pick
sections defined on a finite family of open sets which
cover $M$, and study an analog of the cocycle for those sections.
The key to obtaining the estimate in Equation \ref{usest}
lies in estimating the effects of
conjugation in $H$ by certain elements.  

We begin by introducing some notation.  Suppose $G$ acts by automorphisms
of a principal $H$-bundle $P \to M$, and $V$ is an $H$-space with $E_V
\to M$ the associated bundle, i.e., $E_V = (P \times V)/H$.
Let $C^r(M;E_V)$ be the set of $C^r$ sections of the associated bundle
$E_V$ over $M$.  The following is a well-known result which can be
found in \cite{bw}.

\begin{prop}
\label{toppropbicor}
There exists a natural bijective correspondence between
$H$-equivariant $C^r$ maps $P \rightarrow V$ and $C^r(M;E_V)$,
with respect to which $H$-equivariant $G$-invariant $C^r$ maps
$P \rightarrow V$ correspond to $G$-invariant elements of $C^r(M;E_V)$.
Further, similar results hold if we replace $C^r$ with H\"{o}lder, measurable,
etc.
\end{prop}

By Proposition \ref{toppropbicor}, the assumption that there exists 
a measurable $G$-invariant section of
the bundle $E_V \to M$ is equivalent to the existence of a measurable 
$G$-invariant $H$-equivariant map $\Phi:P \to V$.
Since $V$ is a transitive $H$ space with an 
$H$-invariant metric, we may write $V \cong H_v\backslash H$ where $H_v$ is the
isotropy of $H$ for some fixed $v \in V$.  Note that $H_v$ is compact.
We may assume that $\Phi$ is defined on a $G$ 
and $H$-invariant set. To prove the
theorem, we will show that $\Phi$ can be extended to all of $P$ in a 
H\"{o}lder manner.  The proof consists of two parts.  In the 
first part, we analyze the action of a regular element $a \in G$ on $P$
over a stable manifold of $a$ in $M$.  This leads to the 
crucial estimate, Equation \ref{usest}.
Then, using this estimate, we adapt Liv\v{s}ic's methods to obtain
our result. 

Fix a Riemannian metric $\| \cdot \|_M$ on $M$ and an $H$-invariant
Riemannian metric $\|\cdot\|_P$ on $P$.  
Choose finitely many open neighborhoods $U_j \subset M$ and a compact
set \label{defofK} $K \subset P$ such that
\begin{enumerate}
\item[(a)] $\cup_j U_j = M$,
\item[(b)] there exists $\zeta > 0$ such that the $\zeta$ ball around any
point in $M$ is contained in some $U_j$, and
\item[(c)] there exist smooth sections $s_j:U_j \to K\subset P$ with uniform
Lipschitz constant $\delta$.
\end{enumerate}
For any $y^* \in P$ lying in the fiber over $y \in U_i \subset
M$ there exists a unique $h_i(y^*)$ such that $y^* = s_i(y)h_i(y^*)$.  Note
that $h_i$ varies as smoothly as $y^*$ does, for all $y^*$ lying over $U_i$,
and that the $h_i$'s are uniformly Lipschitz on any set of the form
$s_i(U_i)\cdot K_0$ where $K_0 \subset H$ is a compact subset.

Using the definition of a normally hyperbolic
diffeomorphism in Equation \ref{anosoveqn}, there exists some $k$ such that 
$Ce^{-kA} < 1$.
For ease of notation, we will replace $a^k$ with $a$ and $kA$ with $A$,
i.e., we will assume that $Ce^{-A} < 1.$
 
Let $x \in M$.  Then by our choice of the $U_j$'s, for every $m$ 
there exists some $j$ such that the $\zeta$ ball about $a^mx$ lies entirely
within $U_j$.  For any $x$ and $m$, pick $i(m,x)$ to be such a $j$. 
We remark that although there is no canonical choice for these $i(m,x)$,
the calculations we desire will not depend on these choices.
These $h_{i(m,x)}(a^{m}x^*)$ will serve as our cocycle analogs, and
we want to estimate $h_{i(m,x)}(a^{m}x^*)$ as $x$ varies along
a stable manifold of $a$ by expanding $h_{i(m,x)}(a^{m}x^*)$
as a product where each factor will be controllable.
To this end, we set $q(a,0,x^*) = h_{i(0,x)}(x^*)$ and $q(a,j,x^*) = 
h_{i(j,x)}(a^{j}x^*)[h_{i(j-1,x)}(a^{j-1}x^*)]^{-1}$ for 
$j \ge 1$.
Hence, we can write 
\begin{eqnarray}
\label{qheqn}
h_{i(j,x)}(a^{j}x^*) = q(a,j,x^*)q(a,j-1,x^*) \cdots q(a,0,x^*) 
\end{eqnarray}
for any $j>0$.

With the aim of deriving some preliminary properties of the 
$q(a,j,x^*)$, we use the definition of the $h_i$'s to obtain
\begin{eqnarray*}
a^{j}x^* & = & a[a^{j-1}x^*]= a[s_{i(j-1,x)}(a^{j-1}x)
h_{i(j-1,x)}(a^{j-1}x^*)] \\
& = & [a s_{i(j-1,x)}(a^{j-1}x)]h_{i(j-1,x)}(a^{j-1}x^*) \\
& = & [s_{i(j,x)}(a^{j}x)h_{i(j,x)}(as_{i(j-1,x)}(a^{j-1}x))]
h_{i(j-1,x)}(a^{j-1}x^*).
\end{eqnarray*}
By uniqueness, we conclude
\[ h_{i(j,x)}(a^{j}x^*) = h_{i(j,x)}(as_{i(j-1,x)}(a^{j-1}x))
\cdot h_{i(j-1,x)}(a^{j-1}x^*), \]
and therefore, for any $j >0$,
\begin{eqnarray*}
q(a,j,x^*) & = & h_{i(j,x)}(a^{j}x^*)
[h_{i(j-1,x)}(a^{j-1}x^*)]^{-1} \\
& = &  h_{i(j,x)}(as_{i(j-1,x)}(a^{j-1}x))
h_{i(j-1,x)}(a^{j-1}x^*)[h_{i(j-1,x)}(a^{j-1}x^*)]^{-1} \\
& = & h_{i(j,x)}(as_{i(j-1,x)}(a^{j-1}x)).
\end{eqnarray*}
So, although $q(a,j,x^*)$ is defined as a function of $x^*$,
for $j>0$ it actually depends only on $x$ and our choice of $i(m,x)$.
Of course, $q(a,0,x^*)=
h_{i(0,x)}(x^*)$ still depends on $x^*$.  For $j>0$, let us 
define $q(a,j,x) = q(a,j,x^*)$.  So, for $j>0$, we get 
\begin{eqnarray}
\label{qinveqn}
q(a,j,x)= h_{i(j,x)}(as_{i(j-1,x)}(a^{j-1}x)).
\end{eqnarray}
Since the $s_i$'s take values in a compact set, it follows
that there exists a compact set $K_1\subset H$ such that 
$K_1^{-1} \subset K_1$ and $q(a,j,x) \in K_1$ for
every $j >0$ and for every $x\in M$.  Without loss of generality we
may assume that $\mathcal{K} \times \{0\} \subset K_1$.

Choose $y \in W^s_a(x)$ and let $y^* \in P$ lie in the fiber over $y$.  
Using Equation \ref{anosoveqn}, we get
\[ d_M(a^{j}x,a^{j}y) < C e^{-jA}d_s(x,y), \]
where $d_M$ is the metric in $M$, and $d_s$ is the induced metric on the
leaves of the stable foliation.  In particular, 
if $d_s(x,y) < \zeta$, then by choice of $i(j,x)$, 
$a^{j}y \in U_{i(j,x)}$ for all $j \ge 0$.
Pick an inner product $\left<\cdot,\cdot\right>$ on $\frak{h}$, the Lie algebra
of $H$, which is invariant under the compact group $Ad(H_v)$, 
and let $d_H$ be the corresponding left invariant Riemannian metric on $H$.  

Since the $h_i$'s are uniformly Lipschitz on any set of the form 
$s_i(U_i) \cdot K_0$ where $K_0 \subset H$ is a compact set, and since
$aK_0 \subset H$ can be written as a finite
union of sets of this form, it follows that the $h_i$'s
are uniformly Lipschitz on $aK$ with constant, say $\lambda$, which
depends only on $a$.  Let $c$ and $\theta$ be the H\"{o}lder
constant and exponent for 
multiplying $K$ by $a$ in $P$, so that for any $p_1,p_2 \in K$, we have
$d_P(ap_1,ap_2)  \le cd_P(p_1,p_2)^\theta$.  If we let $\Delta =
\lambda c \delta^\theta C^\theta$, then for all $j>0$,
\begin{equation}
\label{H1}
\begin{array}{c}
d_H(q(a,j,x^*),q(a,j,y^*)) \\[12pt]
\le d_H(h_{i(j,x)}(as_{i(j-1,x)}(a^{j-1}x)),
h_{i(j,x)}(as_{i(j-1,x)}(a^{j-1}y))) \\[12pt]
\le \lambda d_P(as_{i(j-1,x)}(a^{j-1}x),
as_{i(j-1,x)}(a^{j-1}y)) \\[12pt]
\le \lambda c d_P(s_{i(j-1,x)}(a^{j-1}x),
s_{i(j-1,x)}(a^{j-1}y))^\theta \\[12pt]
\le \lambda c \delta^\theta
d_M(a^{j-1}x,a^{j-1}y)^\theta\\[12pt]
\le \lambda c \delta^\theta C^\theta e^{-\theta (j-1)A}
d_s(x,y)^\theta.\\[12pt]
= \Delta e^{-\theta (j-1)A}d_s(x,y)^\theta.
\end{array}
\end{equation}
It is crucial to 
note that although $\Delta$ depends on $a$, it does not depend on $j$.

At this point, we will begin to analyze how $h_{i(n,x)}(a^{m}x^*)$
varies over a stable manifold.
Using Equation \ref{qheqn}, the triangle inequality,
and left invariance we have
\[ d_H(h_{i(m,x)}(a^{m}x^*),h_{i(m,y)}(a^{m}y^*)) \]
\[ = d_H(q(a,m,x^*) \cdots q(a,0,x^*), q(a,m,y^*) \cdots q(a,0,y^*)) \]
\begin{eqnarray*}
\lefteqn{\le  \sum_{j=0}^m d_H(q(a,m,y^*) \cdots q(a,j+1,y^*)q(a,j,x^*)
\cdots q(a,0,x^*),} \hspace{0.3in}\\ 
& & q(a,m,y^*) \cdots q(a,j,y^*)q(a,j-1,x^*)\cdots q(a,0,x^*)) 
\end{eqnarray*}
\[ = \sum_{j=0}^m d_H(q(a,j,x^*) \cdots q(a,0,x^*),
q(a,j,y^*) q(a,j-1,x^*)\cdots q(a,0,x^*)) \]
\[ = \sum_{j=0}^m d_H(C_0 \circ \cdots \circ C_{j-1}(q(a,j,x^*)),
C_0 \circ\cdots\circ C_{j-1}(q(a,j,y^*))) \]
where $C_i$ is conjugation in $H$ by $q(a,i,x^*)$ for $i \ge 0$ and
$C_{-1}$ is the identity transformation.

So, in order to obtain an estimate for 
$d_H(h_{i(m,x)}(a^{m}x^*),h_{i(m,y)}(a^{m}y^*))$, we will produce one
for $d_H(C_0 \circ \cdots \circ C_{j-1}(q(a,j,x^*)),
C_0 \circ\cdots\circ C_{j-1}(q(a,j,y^*)))$.
Recall that by assumption $H = \mathcal{K} \ltimes \mathcal{A}$, where
$\mathcal{K}$ is compact and $\mathcal{A}$ is abelian.  Without loss of
generality we may assume that with respect to the left invariant metric
on $H$, that $\mathcal{K}$ acts via isometries on $\mathcal{A}$.  For
$k_i\in\mathcal{K}$ and $a_i\in \mathcal{A}$, we have that multiplication
is given by $(k_1,a_1)(k_2,a_2)= (k_1k_2,k_2^{-1}a_1+a_2)$.

Let $q(a,i,x^*)=(r_i,s_i)\in\mathcal{K}\ltimes\mathcal{A}$. Then for
$h=(k,a)\in\mathcal{K}\ltimes\mathcal{A}$, we have
\[ C_0 \circ \cdots \circ C_{j-1}((k,a))
 = (r_0,s_0)\cdots(r_{j-1},s_{j-1})(k,a)(r_{j-1},s_{j-1})^{-1}\cdots
(r_0,s_0)^{-1} \]
\[=\Big((r_0\cdots r_{j-1})k(r_0\cdots r_{j-1})^{-1}, r_0
\Big\{ \left[
(r_1\cdots r_{j-1})(k^{-1}-Id)(r_1\cdots r_{j-1})^{-1} \right] s_0\]
\[+r_1\Big\{ \left[ (r_2\cdots r_{j-1})(k^{-1}-Id)(r_2\cdots r_{j-1})^{-1})
\right]s_1 + \cdots + r_{j-1}\left\{ [k^{-1}-Id]s_{j-1} + a \Big\}\cdots
\Big\}\right\}\Big)\]

Since $\mathcal{K}$ is compact, there exists some constant $J>1$ such that
for any $k_1,k_2 \in \mathcal{K}$, $d_H((k_1k_2k_1^{-1},0),(1,0)) \le 
Jd_H((k_2,0),(1,0))$.  Additionally, since $K_1$ is compact, 
there exists some constant $L>1$
such that $d_H((k,a),(1,0))<L$ for every $(k,a)\in K_1$.  In particular, 
$d_H(q(a,j,x^*),(1,0))<L$ for all $j$ and almost every $x\in M$.
Further, there exists some $\Upsilon>1$
so that $d_H((1,(k^{-1}-Id)a),(1,0)) \le 
\Upsilon d_H((k,0),(1,0))d_H((1,a),(1,0))$ for all 
$k\in\mathcal{K}$.  Also, there exists some constant $\omega>1$ such that
$d_H((k,0),(1,0)),d_H((1,a),(1,0)) < \omega d_H((k,a),(1,0))$ for all
$(k,a) \in K_1$.  
Finally, by left invariance
and the triangle inequality note that $d_H((k,a),(1,0)) \le d_H((k,0),(1,0))+ 
d_H(((1,a),(1,0))$.

Using the fact that $\mathcal{K}$ acts via isometries on $\mathcal{A}$,
we have for all $j>0$,
\begin{equation*}
\begin{array}{c}
d_H(C_0 \circ \cdots \circ C_{j-1}((k,a)),(1,0))\\[12pt]
\le d_H(([r_0\cdots r_{j-1}]k[r_0\cdots r_{j-1}]^{-1},0),(1,0)) \\[12pt]
+ \left(\sum_{i=0}^{j-1} d_H(1,r_0\cdots r_i([r_{i+1}\cdots r_{j-1}
(k^{-1}- Id) [r_{i+1}\cdots r_{j-1}]^{-1})s_i,(1,0))\right)\\[12pt]
+ d_H((1,r_0 \cdots r_{j-1}(a)),(1,0))\\[12pt]
\le Jd_H((k,0),(1,0))
+ \left(\sum_{i=0}^{j-1} d_H((1,[r_{i+1}\cdots r_{j-1}](k^{-1} - Id)[r_{i+1}
\cdots r_{j-1}]^{-1}s_i),(1,0))\right)\\[12pt]
 + d_H((1,a),(1,0))\\[12pt]
\le Jd_H((k,0),(1,0))\\[12pt]
+\left(\sum_{i=0}^{j-1} \Upsilon d_H((1,[r_{i+1}\cdots r_{j-1}]k[r_{i+1}
\cdots r_{j-1}]^{-1}s_i),(1,0))d_H((1,s_i),(1,0))\right)\\[12pt] 
+d_H((1,a),(1,0))\\[12pt]
\le Jd_H((k,0),(1,0))
+\left(\sum_{i=0}^{j-1} \Upsilon J d_H((k,0),(1,0))d_H((1,s_i),(1,0))\right)
+ d_H((1,a),(1,0))\\[12pt]
\le Jd_H((k,0),(1,0))
+\left(\sum_{i=0}^{j-1} \Upsilon JLd_H((k,0),(1,0)) \right) 
+d_H((1,a),(1,0)) \\[12pt]
\le (j+2)\Upsilon\omega JL d_H((k,a),(1,0)).
\end{array}
\end{equation*}
By Equation \ref{H1}, 
\[d_H(q(a,j,x^*),q(a,j,y^*))= d_H(q(a,j,y^*)^{-1}
q(a,j,x^*),(1,0)) \le \Delta e^{-\theta (j-1)A}d_s(x,y)^\theta.\]
Hence, we have for almost every $x\in M$ and for $y\in W_a^s(x)$,
\begin{equation}
\label{FF1}
\begin{array}{c}
d_H(h_{i(m,x)}(a^{m}x^*),h_{i(m,y)}(a^{m}y^*))\\[12pt] 
\le \sum_{j=0}^m d_H(C_0 \circ \cdots \circ C_{j-1}(q(a,j,x^*)),
C_0 \circ\cdots\circ C_{j-1}(q(a,j,y^*))) \\[12pt] 
\le \sum_{j=0}^m d_H(C_0 \circ \cdots \circ C_{j-1}(q(a,j,y^*)^{-1}
q(a,j,x^*)),(1,0)) \\[12pt] 
\le d_H(h_{i(0,x)}(x^*),h_{i(0,y)}(y^*)) +
\left(\sum_{j=1}^m (j+2)\Upsilon\omega JL \Delta 
e^{-\theta (j-1)A}\right) d_s(x,y)^\theta.  
\end{array}
\end{equation}
Since $\theta,A>0$, $\sum_{j=1}^\infty (j+2)\Upsilon\omega JL \Delta 
e^{-\theta (j-1)A}$ converges absolutely. 

Summarizing, 
there exists a constant $Q>0$ such that for every $m>0$, 
and every $y \in W^s_\zeta(x)$, the $\zeta$-ball
in $W^s(x)$ about $x$,
\begin{equation}
\begin{array}{c}
\label{stableest}
d_H(h_{i(m,x)}(a^{m}x^*),h_{i(m,x)}(a^{m}y^*)) \\ [12pt]
< d_H(h_{i(0,x)}(x^*),h_{i(0,x)}(y^*)) +
Q d_s(x,y)^\theta. 
\end{array}
\end{equation}
Similarly, 
there exists a constant $Q' >0$ such that for every
$m<0$, 
and every $y \in W^u_\zeta(x)$, the $\zeta$-ball
in $W^u(x)$ about $x$,
\begin{equation}
\begin{array}{c}
\label{unstableest}
d_H(h_{i(m,x)}(a^{m}x^*),h_{i(m,x)}(a^{m}y^*)) \\ [12pt]
< d_H(h_{i(0,x)}(x^*),h_{i(0,x)}(y^*)) + 
Q' d_u(x,y)^\theta. 
\end{array}
\end{equation}
Without loss of generality we can assume that $Q = Q'$.
We may therefore conclude that if $x^*, y^*$ are 
contained in $K \subset P$ (recall that $K$ is defined on 
page~\pageref{defofK}), then
there exists a bound $\kappa$ (depending on $K$) such that
\begin{eqnarray}
\label{usest}
d_H(h_{i(m,x)}(a^{m}x^*),h_{i(m,x)}(a^{m}y^*)) < \kappa d_P(x^*,y^*)^\theta
\end{eqnarray}
provided that $y \in W^s_\zeta(x)$ or $y \in W^u_\zeta(x)$ and
$x \in M_s \cap M_u$.

We return to our measurable $G$-invariant $H$-equivariant
map $\Phi: P \to V$.  By $G$-invariance and $H$-equivariance, we have
\begin{equation}
\label{ghequiv}
\begin{array}{c}
\Phi(x^*) = \Phi(a^{m}x^*)= \Phi(s_{i(m,x)}(a^{m}x)h_{i(m,x)}(a^{j}x))\\
[12pt] = [h_{i(m,x)}(a^{j}x)]^{-1}\Phi(s_{i(m,x)}(a^{m}x)).
\end{array}
\end{equation}
As mentioned above, we need to show that $\Phi$
can be extended to a H\"{o}lder function on all of $P$.  
We will do this by adapting Liv\v{s}ic's argument in
\cite[Theorem 9]{livsic1}.  

Using Lusin's Theorem \cite[Theorem 2.24]{rudin1}, there exists
a compact set $M' \subset M$ with measure $> 1/2$ such that
$M' \cap U_i$ is compact and $\Phi\circ s_i:M' \cap U_i \to H$ 
is uniformly continuous for every $i$.  
By $H$-equivariance, it will suffice to show that $\Phi$ is H\"{o}lder
on $K$, since $K$ is a compact set in $P$ containing an open set
covering all of $M$.

Let $M_0$ be the set in $M$ consisting of all $x$ such that 
\begin{enumerate}
\item[(a)] $\displaystyle \lim_{h\to \infty} \frac{1}{h} \sum_{i=0}^h 
\chi_{M'}(a^ix) = \mu(M') $
\item[(b)] $\displaystyle \lim_{h\to -\infty} \frac{1}{-h} \sum_{i=h}^0 
\chi_{M'}(a^ix) = \mu(M')$, and
\item[(c)] $x$ lies in the projection from $P$ to $M$ of the set on
which $\Phi$ is $G$-invariant and $H$-equivariant.
\end{enumerate}
Let $\mu_s$ and $\mu_u$ denote the conditional measures on the
stable and unstable foliations.  Recursively define
\[ M_{n+1} = \{ x \in M_n | \mu_s(W^s_\zeta(x) \setminus M_n) = 
\mu_u(W^u_\zeta(x) \setminus M_n)=0 \}, \]
and let $M_\infty = \cap_n M_n$.
By absolute continuity of the stable and unstable foliations we have
$\mu(M_\infty) = 1$.  By \cite{livsic1}, absolute continuity also implies
that there exist constants $\xi_1,\omega>0$
such that if $d_M(x_0,x_5) < \omega$ for $x_0,x_5 \in M_\infty$, 
then there exist
points $x_1,x_2,x_3,x_4 \in M_\infty$ such that 
\begin{enumerate}
\item[(a)] $x_1 \in W^s(x_0), d_s(x_0,x_1) < \xi_1 d_M(x_0,x_5),$
\item[(b)] $x_2 \in W^u(x_1), d_u(x_1,x_2) < \xi_1 d_M(x_0,x_5),$
\item[(c)] $x_3 \in W^s(x_2), d_s(x_2,x_3) < \xi_1 d_M(x_0,x_5),$ 
\item[(d)] $x_3 \in W^u(x_4), d_u(x_3,x_4) < \xi_1 d_M(x_0,x_5),$ and
\item[(e)] $x_4$ lies in the $G$ orbit of $x_5$, $d_M(x_4,x_5)< 
\xi_1 d_M(x_0,x_5)$.
\end{enumerate}

To demonstrate that $\Phi$ can be extended to a H\"{o}lder function, we
want to pick $x_0^*$ and $x_5^*$ lying in the fibers over $x_0$ and $x_5$
and show that $d_V(\Phi(x_0^*),\Phi(x_5^*))$ can be bounded by 
$d_H(x_0^*,x_5^*)^\theta$.  This will require the use of Equation \ref{usest},
which is valid only for points lying in the same stable or 
unstable manifolds.  In addition to using canonical coordinates, we 
will require our points to lie in $M_\infty$, which necessitates
the use of the points $x_1,x_2,x_3$, and $x_4$.

For given $x_0^*,x_5^* \in K$ lying in the fibers over $x_0, x_5$, 
it is possible to pick $x_i^* \in K, i=1,2,3$ in 
the fibers over $x_i$ such that 
\[ d_P(x^*_i,x^*_{i+1}) < \xi_2 d_P(x_0^*,x_5^*) \]
for some $\xi_2 > 1$.  Additionally, we may assume without loss of 
generality that $x_4^*$ has been chosen so that it lies in the $G$-orbit
of $x_5^*$.

Since $\Phi$ is $G$-invariant, we have $\Phi(x_4^*)=\Phi(x^*_5)$, so that
$d_V(\Phi(x_4^*),\Phi(x_5^*))=0$.
Pick $d_P(x_0^*,x_5^*) < \frac{\zeta}{\xi_2}$. Then 
\begin{eqnarray}
\label{sumeqn}
d_V(\Phi(x_0^*),\Phi(x_5^*)) \le \sum_{i=0}^3 d_V(\Phi(x_i^*),\Phi(x_{i+1}^*)). 
\end{eqnarray}

Now, using Equation \ref{ghequiv}, the triangle inequality, and left 
invariance, for $i=0,1,2,3$, we have
\begin{equation*}
\begin{array}{c}
 d_V(\Phi(x_i^*),\Phi(x_{i+1}^*)) \\ [12pt]
= d_V([h_{i(m,x_i)}(a^{m}x_i)]^{-1}\Phi(s_{i(m,x_i)}(a^{m}x_i)),
[h_{i(m,x_{i})}(a^{m}x_{i+1})]^{-1}\Phi(s_{i(m,x_{i})}(a^{m}x_{i+1}))\\
[12pt] \le d_V(\Phi(s_{i(m,x_i)}(a^{m}x_i)),\Phi(s_{i(m,x_{i})}
(a^{m}x_{i+1}))) \\ [12pt]
+ d_V(\Phi(s_{i(m,x_{i})}(a^{m}x_{i+1})),
[h_{i(m,x_i)}(a^{m}x_i)][h_{i(m,x_{i})}(a^{m}x_{i+1})]^{-1}
\Phi(s_{i(m,x_{i})}(a^{m}x_{i+1}))).
\end{array}
\end{equation*}

First, note that if $i$ is odd (even) then
as $m \to \infty (-\infty)$,
$a^{m}x_i$ and $a^{m}x_{i+1}$ converge.
Thus, using uniform continuity of $\Phi\circ s_i$ on $M'$, it follows that 
\[ d_V(\Phi(s_{i(m,x_i)}(a^{m}x_i)),\Phi(s_{i(m,x_{i})}
(a^{m}x_{i+1})))\]
can be made as small as we like by choosing $m$ sufficiently large.

Next, using Equation \ref{usest}, we can bound
$d_H(1,[h_{i(m,x_i)}(a^{m}x_i)][h_{i(m,x_{i})}(a^{m}x_{i+1})]^{-1})$
in terms of $d_P(x_0^*,x_5^*)^\theta$ independently of $m>0$ ($m<0$ if $i$ is
even).  Since $\Phi \circ s_i$ takes values in a compact set, it follows
that we can bound
\[ d_V(\Phi(s_{i(m,x_{i})}(a^{m}x_{i+1})),
[h_{i(m,x_i)}(a^{m}x_i)][h_{i(m,x_{i})}(a^{m}x_{i+1})]^{-1}
\Phi(s_{i(m,x_{i})}(a^{m}x_{i+1}))) \]
in terms of $d_P(x_0^*,x_5^*)^\theta$ independently of $m>0$ ($m<0$ if $i$ is
even).

Hence, for each $i$, we have 
a bound for $d_H(\Phi(x_i^*),\Phi(x_{i+1}^*))$ in terms of 
$d_P(x_0^*,x_5^*)^\theta$.
Consequently, there exists a bound for $d_V(\Phi(x_0^*),\Phi(x_5^*))$
in terms of $d_P(x_0^*,x_5^*)^\theta$.
That $\Phi$ can be H\"{o}lder extended to all of $K$ follows using this
bound, which completes the proof.
\qed

\section{Anosov Actions of Semisimple Groups}

Our goal in this section is to combine Theorem \ref{livthm} 
with Topological Superrigidity to obtain geometric
information for Anosov actions of semisimple Lie groups, 
and, in particular, Theorem \ref{secthm}.
We begin this section by stating these results as well as
briefly describing the relevant definitions and facts necessary
to understand them.  See \cite{erg1} and \cite{z9} for a 
complete discussion of the related concepts and properties.
We then continue in the following subsections
by proving our results.

\subsection{Statement of Results}

Suppose $G$ is a connected semisimple Lie group without
compact factors such that each simple factor of $G$ has
$\Bbb R$-rank at least 2.   Suppose that $G$ acts on a closed manifold
$M$ such that the $G$ action is Anosov and volume preserving.
Then $G$ acts via derivatives on the general frame bundle over $M$, and
since the action preserves a volume, there exists a $G$-invariant reduction
to a principal $H$ bundle with $H \subset SL(n,\Bbb R)$, an
$\Bbb R$-algebraic subgroup.  It follows therefore that any algebraic 
hull of the $G$ action on the general frame bundle is contained in 
$SL(n,\Bbb R)$.  In fact, we have the following very precise description
of the algebraic hull.

\begin{thm}\label{THM2}
Suppose $G$ is a connected semisimple Lie group without
compact factors such that each simple factor of $G$ has
$\Bbb R$-rank at least 2.   Suppose that $G$ acts on a closed manifold
$M$ such that the $G$ action is Anosov and volume preserving.
Then the H\"{o}lder algebraic hull of the $G$ action by
derivatives on the frame bundle is reductive with
compact center.
\end{thm}
\begin{remark}
This theorem also holds for volume preserving Anosov actions of a
cocompact lattice $\Gamma \subset G$ where each simple factor of $G$
has rank 2 or 3.  Additionally, this result holds whenever there are
enough normally hyperbolic elements in $\Gamma$.  For instance, if there
exists a split Cartan $A\subset G$ such that $A \cap \Gamma$ is cocompact
in $A$ and such that  $\{ \log{\gamma} \| \gamma \in A \cap \Gamma \}$
is dense in the space of directions of $\frak{a}$, the Lie algebra for
$A$.  Or, if $\Gamma$ is invariant under the Weyl group for $A$.  The main
issue is to ensure that an appropriate version of
Proposition \ref{rootprop} below holds.
\end{remark}

Let $G$ act on a set $S$.  Then $s \in S$
is called a {\em parabolic invariant} if there exists a parabolic
subgroup $Q \subset G$ such that $Q$ fixes $s$.  In particular, if
$G$ acts by automorphisms of a principal $H$-bundle $P \to M$, and $V$ is
an $H$-space with $E_V \to M$ the associated bundle, a {\em parabolic
invariant section} is a section of $E_V \to M$ invariant under a parabolic
subgroup of $G$.

Frequently, it is possible to obtain a great deal of
geometric information on a generic subset of $M$.  Of course, we also
wish to emphasize the distinction between proper generic subsets and
all of $M$.
Let $G$ act via principal bundle automorphisms on $P(M,H)$ with
H\"{o}lder (measurable) algebraic hull $L$.  The $G$ action on $P(M,H)$ is
{\em H\"{o}lder (measurably) complete} if there exists a
H\"{o}lder (measurable) $G$-invariant section of $E_{H/L}$.
Given a principal $H$-bundle $P \to M$, and $V$ an $H$-space, a section
$\phi$ of $E_V$ is called {\em effective for $P$} if $H$ acts effectively
on $\Phi(P)$ where $\Phi:P \to V$ is the $H$-map corresponding to $\phi$.
Suppose $G$ acts ergodically and H\"{o}lder (measurably) completely on 
$P$ where $H$ is an algebraic group and $V$ is an algebraic variety.
Then a H\"{o}lder (measurable) section $\phi:M \to E_V$ is 
{\em H\"{o}lder (measurably) $G$-effective} 
if it is effective for $P_1 \subset P$ where $P_1$ is a $G$-invariant reduction
to $L \subset H$, the H\"{o}lder (measurable) algebraic hull.

Any $\Bbb R$-split Cartan $A\subset G$ contains a normally hyperbolic
element $a\in A$ (cf. Lemma \ref{nhyperinA}).
Suppose dim$(E^-_a)=k$ and let
$Gr(k,n)$ be the set of $k$-dimensional subspaces of $\Bbb R^n$.
Form the associated
bundle $E_{Gr(k,n)} = (P \times Gr(k,n))/SL(n,\Bbb R)$.  
Then we can define an 
$A$-invariant H\"{o}lder section $\phi:M \to E_{Gr(k,n)}$ by 
setting $\phi(m)= E^-_a(m)$.  The following demonstrates that this section 
actually possesses greater geometric properties.

\begin{prop} \label{THM1}
Suppose $G$ is a connected semisimple Lie group without
compact factors such that each simple factor of $G$ has
$\Bbb R$-rank at least 2.   Suppose $G$ acts on a closed manifold
$M$ such that the $G$ action is Anosov and volume preserving.
Let $H$ be the H\"{o}lder algebraic hull.
Then, modulo a compact subgroup of $H$, there exists a 
H\"{o}lder $G$-effective parabolic invariant section 
$\phi:M \to E_{Gr(k,n)}$, i.e., there exist 
\begin{enumerate}
\item a normal subgroup $N \subset H$ which fixes $\Phi(P/N)$,
\item a parabolic subgroup $Q \subset G$, and 
\item a $Q$-invariant $N\backslash H$-equivariant map $\Phi:P/N
\to Gr(k,n)$,
\end{enumerate}
such that $H/N$ acts effectively on $\Psi(P/N)$.
\end{prop}

Let $P \to M$ be a principal $H$-bundle on which $G$ acts via bundle
automorphisms.  If $\pi: G \to H$ is a homomorphism, then a section
$s:M \to P$ is called {\em totally $\pi$-simple} if for $g \in G$ and
$m \in M$,
\[ s(gm) = g.s(m).\pi(g)^{-1}. \]
Here, of course, $M$ and $P$ are left $G$-spaces, and $P$ is a right 
$H$-space.

\begin{thm}[Topological Superrigidity]
Let $G$ be a connected semisimple Lie group, $\Bbb R$-rank $(G) \ge 2$, with
$G$ acting via bundle automorphisms on $P(M,H)$, $H$ an algebraic
$\Bbb R$-group, $V$ an $\Bbb R$-variety on which $H$ acts algebraically,
and $G$ acting ergodically on $M$ with respect to a probability measure
$\mu$ where supp$(\mu)= M$.  Assume the action is H\"{o}lder complete.  If
there exists a $G$ effective H\"{o}lder parabolic invariant 
section $\psi$ of $E_V \to M$, then, by possibly passing to a finite cover
of $G$, there exist
\begin{enumerate}
\item a homomorphism $\pi :G \to H$,
\item $v_0 \in V$, and
\item a totally $\pi$-simple H\"{o}lder section $s$ of $P \to M$
\end{enumerate}
such that $\psi$ is the associated section $(s,v_0)$ of $E_V$, i.e.
$\psi(x)=[s(x),v_0]$.
\label{TSR}
\end{thm}
\begin{remark}
\label{TSRremark}
\begin{enumerate}
\item Topological Superrigidity was proved by Zimmer in \cite{z9}, and the 
proof of a generalization appears in \cite{erg1}.  The version above differs
from the original version of Topological Superrigidity in that we have
replaced H\"{o}lder functions for $C^r$ functions.  However, the proof for
Theorem \ref{TSR} requires only a minor modification of the proofs presented
in \cite{erg1,z9}.
\item Effectiveness of $\psi$ ensures that we can see enough of $H$ in the
image $\psi$.  Note that by passing to a suitable subquotient, it 
is always possible to obtain an effective section, and, 
in fact, this procedure is required for most of our applications.  
More explicitly, let $\Psi:P \to V$ be the $H$-equivariant map corresponding
to $\psi$, and let $N\subset H$ be the 
kernel of the $H$ action on $\Psi(P)$.
\label{TSRremark2}
We can then obtain a $G$-effective section
of the $V$ associated bundle to the principal bundle $P/N$, and by applying
the last theorem to this section, we obtain a homomorphism $\sigma:G \to H/N$
and a totally $\sigma$-simple H\"{o}lder section $s$ of $P/N \to M$.  Since
$H$ is reductive and $N$ is normal, $H$ is an almost direct product
$H= N\cdot H_1$ of $N$ with a normal subgroup $H_1\subset H$.  Then,
we can produce a homomorphism $\pi:G \to H_1\subset H$ which
modulo $N$ factors to $\sigma$.  Then the section 
$s$ will also be totally $\pi$-simple.
See \cite{erg1} for further details.
\end{enumerate}
\end{remark}

Our main result, which follows, is obtained by combining Proposition \ref{THM1}
and Theorem \ref{TSR}.

\begin{thm}
\label{secthm}
Suppose $G$ is a connected semisimple Lie group of higher rank without
compact factors such that each simple factor of $G$ has
$\Bbb R$-rank at least 2.   Suppose that $G$ acts on a closed manifold
$M$ such that the $G$ action is Anosov and volume preserving.
Let $H$ be the H\"{o}lder algebraic hull of the $G$ action on $P \to M$, 
the $G$-invariant reduction of the derivative action on the full frame
bundle over $M$.
Then, by possibly passing to a finite cover of $G$, there exist
\begin{enumerate}
\item a normal subgroup $K \subset H$,
\item a H\"{o}lder section $s:M \to P/K$, and
\item a homomorphism $\pi:G \to H$ (obtained from a homomorphism $G\to H/K$
as in Remark \ref{TSRremark}.\ref{TSRremark2}),
\end{enumerate}
such that $s(gm) = g.s(m).\pi(g)^{-1}$ for every $g \in G$ and every $m \in M$.

Moreover, if $\pi$ is multiplicity free, i.e., all irreducible 
subrepresentations of $\pi$ have multiplicity one,
then $K\subset H$ is a compact normal subgroup.
\end{thm}

\begin{cor}
\label{holdermetriccor}
Let $G$, $P$, $M$, and $H$ be as in Theorem \ref{secthm}.  Assume the
irreducible subrepresentations of $\pi$ are multiplicity free so that
$K$ is compact.
Let $A$ be a maximal $\Bbb R$-split Cartan of $G$, with $\{ \chi \}$
the set of weights of $\pi$ with respect to $A$.
There exist 
\begin{enumerate}
\item a $K$-invariant H\"{o}lder Riemannian metric, $\|\cdot\|_K$, 
on $M$, and 
\item a $K$-invariant H\"{o}lder decomposition $TM = \bigoplus E_\chi$
\end{enumerate}
such that for every $v \in E_\chi$ and $a \in A$,
\[ \|av\|_K= e^{\chi(\log{a})}\|v\|_K. \]
\end{cor}
\begin{proof}
Every frame $f$ over a point $p \in M$ determines an inner
product $\langle\cdot,\cdot \rangle_f$ on $T_p M$ by declaring $f$ to
be orthonormal. Then a $K$-orbit of frames $fK$ determines an
inner product $\langle \cdot,\cdot \rangle_{fK}$ by averaging the
$\langle\cdot,\cdot\rangle_{fk}$ over $K$. This defines a smooth
map from $P/K$ to the bundle of $K$-invariant inner products on
tangent spaces over $M$. Composing this map with the H\"{o}lder
section $s:M \rightarrow P/K$ determines a $K$-invariant H\"{o}lder
Riemannian metric $\parallel \cdot \parallel _K$ on $M$. Pick a
measurable section $\sigma : M \rightarrow P$ which projects onto $s$.
Then $\sigma$ is a measurable totally $\pi$-simple framing of $M$.
Let $TM =\oplus E_i$ be the Lyapunov decomposition corresponding to
the action of $A$. The exponents for the Lyapunov decomposition
correspond to the weights of $A$ under $\pi$
(cf. e.g. \cite{feres,furstenberg,erg2}).
Thus we can write $TM = \bigoplus E_\chi$
such that for every $v \in E_\chi$ and $a \in A$,
$\|av\|_K= e^{\chi(\log{a})}\|v\|_K$.

It remains to see that the $E_\chi$ are H\"{o}lder.  However, 
the Lyapunov decomposition is determined by the measurable section
$\sigma$ and the weight
spaces for $\pi$ with respect to $A$,  and since $\pi(G)$ and
$K$ commute, it follows that the Lyapunov decomposition
is $K$-invariant, and can therefore be defined by the H\"{o}lder
section $s$ instead.  It follows that the $E_\chi$ are H\"{o}lder.
\end{proof}

\subsection{Proof of Theorem \ref{THM2}}
The outline of our proof follows the same basic steps as Zimmer's proof for
the measurable algebraic hull in \cite{z10}.  In the course of our proof,
however, we frequently encounter situations where the gap between a
measurable statement and its H\"{o}lder equivalent become significant.
For example, in the proof that the unipotent radical is trivial, Zimmer
can use Oseledec's Theorem to see that a certain cocycle is exponential, but
there seems to be no H\"{o}lder analog for this.
In these situations, we employ Theorem \ref{livthm} as well as 
finite dimensional representation theory to bridge these gaps.

With the following lemma, we may assume the H\"{o}lder hull 
is Zariski connected. 
\begin{lemma}
Suppose $G$ acts on a principal $H$ bundle $R \to M$ and that
the H\"{o}lder algebraic hull of the $G$ action on $R$ is $H$.  Then
there exists a finite cover $M' \to M$ such that the H\"{o}lder algebraic
hull of the $G$ action on $R \to M'$ is $H^0$, the Zariski connected component
of $H$.
\end{lemma}
\begin{proof}
Let $\displaystyle M' = \frac{(R \times H^0\backslash H)}{H}$.  
Then $M' \to M$ is a finite cover, and $R \to M'$ is a principal $H^0$
bundle on which $G$ acts via bundle automorphisms.
Let $L$ be the H\"{o}lder algebraic hull of the $G$ 
action on $R \to M'$. Clearly $L \subset H^0$, so it remains to show
$H^0 \subset L$.

By the definition of algebraic hull, there exists 
a H\"{o}lder $G$-invariant $H^0$-equivariant 
map $\Phi: R \to L\backslash H^0$.  
Define 
\[ \tilde{\Phi}: \frac{R \times H}{H^0} \to 
\frac{L \backslash H^0 \times H}{H^0} \cong L\backslash H \]
so that $\tilde{\Phi}([p,h]) = [\Phi(p),h]$.  Here, of course, $R$ and
$L \backslash H^0$ are right $H^0$ spaces, and $H$ is a left $H^0$ space.
The equivalence of $H$ spaces,
$\displaystyle\frac{L \backslash H^0 \times H}{H^0} \cong L\backslash H$,
is obtained from the obvious bijective $H$-equivariant map.
With right multiplication by $H$ on itself, both 
$\displaystyle\frac{R \times H}{H^0}$ and 
$\displaystyle\frac{L \backslash H^0 \times H}{H^0}$ are right $H$ spaces, and
$\tilde{\Phi}$ is an $H$-equivariant map.

Let $E=\text{Functions}(H^0 \backslash H, L \backslash H)$,
the set of {\em all} functions from the finite set
$H^0\backslash H$ to $L \backslash H$, and define
$\Psi: R \to E$ by setting $\Psi (p)[h] = \tilde{\Phi}([ph^{-1},h])$.
If $h_0 \in H^0$ then
\[ \Psi(p)[h_0h]=\tilde{\Phi}([ph^{-1}h_0^{-1},h_0h]) = 
\tilde{\Phi}([ph^{-1},h]) = \Psi(p)[h] \]
so that $\Psi(p)$ is in fact well defined.  If we define the $H$ action on
$E$ so that, for $f \in E$ and $k \in H$, $(f.k)([h])= f([hk^{-1}])k$, then
$\Psi$ is also $H$-equivariant:
\begin{equation*}
\begin{array}{c}
(\Psi(p).k)([h]) = \Psi(p)([hk^{-1}])k = \tilde{\Phi}([pkh^{-1},hk^{-1}])k \\
[12pt] = \tilde{\Phi}([pkh^{-1},h]) = \Psi(pk)([h]). 
\end{array}
\end{equation*}
The rest of the argument follows exactly as in the proof of 
\cite[Theorem 9.2.6]{z1}.  Namely, by ergodicity of $G$ on $M$ and
the tameness of the $H$ action on $E$, the image of $\Psi$ is
contained in a single $H$ orbit. We can therefore view $\Psi:R \to H_\phi 
\backslash H$ as a H\"{o}lder $G$-invariant $H$-equivariant map, where
$H_\phi$ is the stabilizer for some $\phi \in E$, which, as the finite
intersection of algebraic subgroups, is algebraic.  Since $H$
is the H\"{o}lder algebraic hull of $R \to M$ and $H_\phi$ is algebraic, we 
must have $H_\phi = H$.  However, as a stabilizer, $H_\phi$ leaves a 
finite subset of 
$L \backslash H$ invariant, thus $L \backslash H$ must itself be finite.
Hence, $H^0 \subset L$.
\end{proof}

For the remainder of the proof we fix some $\Bbb R$-split Cartan
$A\subset G$.  Let $\mathcal{R}$ be the system of roots, and $\mathcal{W}$ the
Weyl group for $A$.
\begin{lemma}  
There exists a normally element $a \in A$ such that $(w(\log{a}),\beta)\ne 0$
for every $w\in \mathcal{W}$ and for every $\beta\in \mathcal{R}$.
\label{nhyperinA}
\end{lemma}
\begin{proof}
By structural stability, the set of normally hyperbolic elements in $G$ is
open.  Therefore, there exists a normally hyperbolic element $g \in G$
which is semisimple.  Let $g=ks$ be the polar decomposition of $g$ with
polar part $s$.  Then $s$ is normally hyperbolic since $k$ is contained in
a compact group and commutes with $s$.  Since some conjugate of $g$ lies
in $A$, it follows $A$ contains a normally hyperbolic element.  By
structural stability, the set of normally hyperbolic elements in $A$ is
open, hence we can find $a\in A$ that satisfies the condition above.
\end{proof}

Henceforth, we work with a fixed element $a \in A$ satisfying
Lemma \ref{nhyperinA}.  Since the volume preserving Anosov actions are
ergodic, the $A$ action is ergodic by Moore's Ergodicity Theorem \cite{z1}.
Since the H\"{o}lder algebraic hull is Zariski connected we can 
write $H = L \ltimes U$ with
$L$ reductive and $U$ unipotent.  Let $P \to M$ be 
the $G$-invariant reduction of the full frame bundle to an $H$-bundle
over $M$.  First, we show that $Z(L)$, the center of $L$,
is compact.

Let $N = [L,L] \ltimes U$.  Then $H/N$ is abelian, so
dividing out by the maximal compact subgroup $C\subset H/N$ 
yields $C \backslash H/N$
which is abelian and contains no compact subgroups.  
Since each simple factor of $G$ has higher rank, $G$ has Kazhdan's property.
By \cite{z1}, the measurable hull of the $G$ action on $(P/N)/C$
must be trivial, so
by Theorem \ref{livthm}, the H\"{o}lder hull is also trivial.  
By the following lemma,
this is a contradiction unless $C = H/N$. It follows, therefore, that
$Z(L)$ must be compact.

\begin{lem}\label{hulllem}
Let $N \subset H$ be a normal subgroup.  If the H\"{o}lder hull of the $G$
action on $P \to M$ is $H$, then the H\"{o}lder hull of the $G$ action on
$P/N \to M$ is $H/N$.
\end{lem}
\begin{proof}
The H\"{o}lder hull for $P/N$ must be contained in $H/N$, but
if there exists a $G$-invariant
H\"{o}lder section
\[ M \to \frac{P/N \times \frac{H/N}{B/N}}{H/N} \cong \frac{P \times H/B}{H},\]
where $B$ is a proper algebraic subgroup of $H$ containing $N$,
then $H$ cannot be the H\"{o}lder hull for $P$.
\end{proof}

It remains, therefore, only to show that $U$ is trivial.  
The first step in this direction will be to show that $U$ must be
contained in the stabilizer for a particular $H$ action.
Since the measurable algebraic hull is contained in the H\"{o}lder algebraic
hull, by \cite{z10}, there exists a measurable framing of $M$ such that
the corresponding derivative cocycle has the form $\kappa(g,m)\pi(g)$ where
$\pi: G \to S$ is a homomorphism and $\kappa:G \times M \to K'$ is a cocycle
taking values in a compact subgroup of $H$ commuting with $\pi(G)$.

Note that $L$ contains the image of the measurable superrigidity 
homomorphism.  Since the action is Anosov and $\pi$ determines the
Lyapunov splitting, the image $\pi(G)$ is a noncompact
semisimple Lie group.  Hence, $L$ is not compact.

Let $L=ZS$
where $S$ is semisimple with no compact factors, $Z$ is compact and 
centralizes  $S$, and the product is almost direct.  If $\frak{U}$ is the
Lie algebra of $U$, then we have a representation
$\rho: L \to GL(\frak{U})$ obtained from the semidirect product.
We will denote the restriction of $\rho$ to $S$ simply as $\rho$.
Also let $\sigma:S \to GL(n,\Bbb R)$ be the representation 
determined by the embedding $S \subset H \subset SL(n,\Bbb R)$.

\begin{lem}
If $U \ne 0$, then $\rho|_S$ is not trivial.
\end{lem}
\begin{proof}
If $\rho|_S$ is trivial, then 
we obtain an amenable algebraic group
$H/(\text{ker}(\rho) \ltimes [U,U]) \cong Z_1 \ltimes U/[U,U]$,
where $Z_1$ is some compact quotient of $Z$.  
By Kazhdan's property, the measurable hull of 
$P/(\text{ker}(\rho) \ltimes [U,U])$ must be contained
in the compact part of this amenable group \cite{z1}.  If we 
could show that the H\"{o}lder hull must also
be contained in the compact part, Lemma \ref{hulllem} would then
force $U/[U,U] = 0$, so that $U=0$. 
As $Z_1$ is compact, there is a $Z_1 \ltimes U/[U,U]$-invariant metric
on $V= (Z_1 \ltimes U/[U,U])/Z_1$.  Hence Theorem \ref{livthm} applied
to this $V$ establishes this, thereby completing the proof.
\end{proof}

Let $\{w_1,\dots,w_t\}$ be some fixed ordering of all the 
elements in $\mathcal{W}$, 
and let $k_i$ be the dimension of $E^-_{w_i(a)}(m)$
and let $l_i$ be the dimension of $E^+_{w_i(a)}(m)$.  Let 
\[V = Gr(k_1,n)\times Gr(l_1,n)\times Gr(k_2,n)\times 
Gr(l_2,n)\times \dots\times Gr(l_t,n),\] 
and define a H\"{o}lder $A$-invariant section $\omega:M \to 
E_V$ so that 
\[ \omega(m) = (E^-_{w_1(a)}(m),E^+_{w_1(a)}(m),\dots,
E^-_{w_t(a)}(m),E^+_{w_t(a)}(m)).\]
Then $\omega$ corresponds to a H\"{o}lder $A$-invariant $H$-equivariant map
$\Omega:P \to V$.  Dividing out by the $H$-action, we
obtain an $A$-invariant map $\tilde{\Omega}:P/H \cong M \to V/H$.
Note that 
$H$ acts tamely on $V$ \cite{z1}.  By ergodicity of $A$ on $M$, 
$\tilde{\Omega}$
is constant, and therefore $\Omega$ is contained in an $H$-orbit.
Thus, we may consider $\Omega$ as a map from $P$ into
$H_{x_0}\backslash H$ where $H_{x_0}$ is the 
stabilizer in $H$ of a point $x_0 \in V$.  

We want to use $\Omega$ to produce a $G$-invariant
$H$-equivariant map, which, in turn, will allow us to calculate the
algebraic hull.  Let $Z$ be the center of $G$ and note that $\Omega$ is
$Z$-invariant. Thus, we can define
$\Omega':P \times (G/ZA) \to H_{x_0}\backslash H$ 
so that $\Omega'(p,g) = \Omega(g^{-1}p)$.
Letting $G$ act on $P \times (G/ZA)$ via the diagonal action, $\Omega'$ becomes
$G$-invariant.  Let $F(G/ZA, H_{x_0}\backslash H)$ be the space of measurable
functions with the topology of convergence in measure, and define
\[ \Psi:P \to F(G/ZA, H_{x_0}\backslash H) \]
by $\Psi(p)(g) = \Omega'(p,g)= \Omega(g^{-1}p)$.

Following the proof of Step 2 in the proof of measurable superrigidity in
\cite{z1} and the proof of Lemma 3.3 in \cite{z10}, we see that $\Psi(p)$ is
a rational function for almost every $p$, i.e., $\Psi: P \to R =
\text{Rat}(G/ZA,H_{x_0}\backslash H)$.
Let $G$ and $H$ act on $R$ so that
\[ (g_1.r)(g) = r(g_1^{-1}g), \text{ and}\]
\[ (r.h)(g) = r(g).h. \]
Then $\Psi$ is $G$ and $H$ equivariant.
By $H$-equivariance of $\Psi$, the degree of $\Psi(p)$ as a rational
function depends only on the base point of $p$ in $M$.  Since $G$ acts
ergodically on $M$, it follows that these rational functions must
have the same degree almost everywhere.
But, the set of rational functions
of a fixed degree is closed, so by continuity we may conclude that
$\Psi(p)$ is rational for every $p$.

Again, following the proof of Step 3 for measurable superrigidity in \cite{z1}
and the comments in \cite{z10}
we can conclude that $\Psi$ is $G$-invariant and contained in a single
$H$ orbit in $R$, i.e.,
$\Psi : P \to H_r\backslash H$ where $H_r$ is 
the stabilizer in $H$ for some $r \in R$.
Since $H_r$ consists of the subgroup of $H$ pointwise fixing the image
of $r$, and since without loss of generality we may assume that 
$x_0 \in V$ lies in the image of
$r$, we have that $H_r \subset H_{x_0}$.  However, by the definition
of algebraic hull
we must have $U \subset H_r$, and consequently $U \subset H_{x_0}$.

To complete the proof of Theorem \ref{THM2}, we will use 
$U \subset H_{x_0}$ to show that $U$ must be trivial.
Since $U\subset H \subset SL(n,\Bbb R)$, we have a natural
action of $\frak{U}$ on $\Bbb R^n$.
Set \[ V_0 = \{ v \in \Bbb R^n |Xv=0 \text{ for all } X \in \frak{U} \}. \]
Since $L$ normalizes $U$, $L$ leaves $V_0$ invariant.  Since $L$ is reductive
with compact center there is a splitting $\Bbb R^n = V_0 \bigoplus V_1$ into
$L$-invariant subspaces.
Since any subgroup of $GL(n)$ fixing all of $\Bbb R^n$ must be
trivial, if we can show that $V_0=\Bbb R^n$, then $U$ must
be trivial.  Hence, we must show that $V_1 = 0$.
To do this, we will show that all of $\frak{U}$ kills the maximal weight
space for any nontrivial irreducible subrepresentation of $\sigma\circ\pi$ 
on $V_1$.
This implies that the given nontrivial irreducible 
subrepresentation cannot lie in
$V_1$, and that $\sigma\circ\pi|_{V_1}$ is a sum of trivial representations.
By the next lemma, the proof of Theorem \ref{THM2}
is complete once we see that 
$\sigma\circ\pi|_{V_1}$ is a sum of trivial representations.
\begin{lemma}
If $\sigma\circ\pi|_{V_1}$ is a sum of trivial representations, then $V_1=0$.
\end{lemma}
\begin{proof}
For $v \in V_1 \setminus \{0\}$, the framing determines a measurable vector
field $\mathcal{V}$ on $M$.  Since $S$ fixes $v$, we have $dg_m(\mathcal{V}(m))
= \kappa(g,m)\mathcal{V}(gm)$.  In particular, with $g= a^n$ and using the
recurrence properties of $a$, it follows that $\mathcal{V}$ has 0 Lyapunov
exponent.  Consequently, $\mathcal{V}$ must lie in the $M_0A$ tangent direction,
where $M_0$ is the compact part of the centralizer of $A$ in $G$
(cf. the beginning of Section 2.1).
Note that $K'$ leaves the 0 Lyapunov direction invariant and therefore the
$M_0A$ tangent direction invariant.  

Fix a basis consisting of left invariant vector fields on $\frak{g}$.  They
determine a framing of the tangent directions to the $G$ orbit on $M$
which under $G$ transforms according to the adjoint representation (cf. the
beginning of Section 2).  Since
the adjoint orbit of any nonzero element in $\frak{m}_0\oplus\frak{a}$
cannot be contained in $\frak{m}_0\oplus\frak{a}$,
there exists some $g \in G$ such that $dg(\mathcal{V}(m))$ cannot
lie in the $M_0A$ direction.  This contradiction establishes that there
cannot be nonzero elements in $V_1$.
\end{proof}

To show that $\sigma\circ\pi|_{V_1}$ is a sum of trivial representations,
we shall need the following three results.
\begin{lem} \label{remA}
Suppose $X \in \frak{U}$ has $\rho\circ\pi$ weight $\nu$, and $v \in V_1$ has 
$\sigma\circ\pi$ weight $\chi$.  Then either $X.v=0$ or $v$
has $\sigma\circ\pi$ weight $\nu + \chi$.
\end{lem}
\begin{proof} 
Let $b$ be in $\frak{a}$, the Lie algebra of $A$.
Then we have
$\pi(b)(X.v) = (X\pi(b) + [\pi(b),X]).v = X(\pi(b).v) + [\pi(b),X].v =
X(\chi\circ\pi(b).v) + (\nu\circ\pi(b))X.v = (\chi\circ\pi(b) + 
\nu\circ\pi(b))X.v$.
\end{proof}

It is important to note that $X.v$ may well lie in a
different irreducible subrepresentation than $v$.
\begin{lem} \label{remB}
Let $X \in \frak{U}$, let $K_a$ 
be the hyperplane $\{ \alpha \in \frak{a}^* | \alpha(\log{a}) = 0 \}$,
let $K_a^+ = \{ \beta \in \frak{a}^* |(\alpha,\log{a}) > 0\}$,
and $K^-_a= \{ \beta \in \frak{a}^* |(\alpha,\log{a}) < 0\}$.
If $v\in \Bbb R^n$ is a weight vector with weight in $K_{w(a)}^\pm$ for any
$w\in \mathcal{W}$, then either 
$X v \in K_{w(a)}^\pm$ or $X v = 0$.
\end{lem}
\begin{proof}
Note that $x_0 \in V$ corresponds to a $2t$-tuple
of linear subspaces 
\[(E^-_{w_1(a)},E^+_{w_1(a)},\dots,E^-_{w_t(a)},E^+_{w_t(a)}).\]
Since $U \subset H_{x_0}$,
$X \in \frak{U}$ implies that $X$ preserves $x_0$, and therefore the
stable and unstable directions for all $w_i(a)$ as well.
\end{proof}

\begin{lem} \label{remC}
Let $N$ be a root for $\frak{S}$, the Lie algebra of $S$, and suppose
$[X, N]=0$ for some $X \in \frak{U}$.  If $X v = 0$,
then $X (Nv) = 0$
\end{lem}
\begin{proof}
This follows since $X(Nv) = [X,N]v + N(Xv) = 0.v + N(0) = 0$.
\end{proof}

Consider an irreducible subrepresentation of $\rho:S \to GL(\frak{U})$ with
maximal weight $\mu$.  Then for any weight vector $X$ in this irreducible, 
its weight has the form $\mu - \sum_i l_i\alpha_i$ 
where the $\alpha_i$ are the 
positive simple roots and the $l_i$ are nonnegative integers.  
This is well known in the complex case \cite{hump1}.  For the real case,
simply consider the complexification.  The only problem is caused by imaginary
root vectors $R_\beta$ which show up in the real case.  However, we
may commute any imaginary $R_\beta$ past all nonimaginary root vectors
$R_\alpha$ since $[R_\beta,R_\alpha]$ is nonimaginary.  Since the maximal
weight space $W_\mu$ of $A$ is invariant under the imaginary root vectors,
the claim follows in the real case.
Moreover $X$ can be written in the form
\[  X = [R_{-\alpha_{i_n}}[ \cdots [R_{-\alpha_{i_1}},X_\mu]\cdots], \]
for some weight vector $X_\mu$ with weight $\mu$ and the $R_{\alpha_i}$ 
nonimaginary root vectors. 
Expanding out the brackets, we obtain
\[ X = \sum_i N_i X_\mu M_i, \] where $N_i$ and $M_i$ have the form
$R_{-\alpha_{i_1}} \cdots R_{-\alpha_{i_n}}$ for varying choices of
$\{\alpha_i\}$.  So to see that all of $\frak{U}$ kills the maximal
weight space for any nontrivial subrepresentation of $\sigma$ on $V_1$,
it clearly suffices to show that $X_\mu$ kills all the weight spaces in 
such an irreducible subrepresentation of $\sigma$ on $V_1$.

Let $\mu$ be the maximal weight for an irreducible subrepresentation
of $\rho$, and let $\alpha$ be the maximal root for $S$.
In particular, we have that $\left<\mu,\alpha\right> > 0$.
Let $K_a,K^+_a,K^-_a$ be as in Lemma \ref{remB}.
By renaming if necessary, we may
assume that the normal vector to $K_a$ lies in the positive Weyl chamber.

\begin{lem}  
Suppose $\beta$ is a weight for $\sigma$ such that $\beta \in K_a^+$
and $w_\alpha(\beta) \in K_a^-$.  
Let $w_\alpha \in \mathcal{W}$ be defined as the reflection through the 
hyperplane perpendicular to $\alpha$.  Then there exists a weight
$\lambda$ in the $\alpha$ string from $w_\alpha(\beta)$ to $\beta$ such that
either
\begin{enumerate}
\item $\lambda \in K^-_a$ and $\lambda + \mu \in K^+_a$, or,
\item $\lambda +\alpha\in K^+_a$ and 
$\lambda + \alpha + w_\alpha(\mu) \in K^-_a$. 
\end{enumerate}
\label{killweightlem}
Further, this result holds if we replace $K_a$ with $K_{w(a)}$ for any
$w\in\mathcal{W}$.
\end{lem}

\begin{proof}
Since $\beta \in K_a^+$ and $w_\alpha(\beta) \in K_a^-$, the 
$\alpha$-string from $\beta$ to $w_\alpha(\beta)$ intersects $K_a$,
say between the weights $\lambda$ and $\lambda + \alpha$.
Form the triangle with vertices $\lambda, \lambda + \mu$, and $\lambda +
\mu - w_\alpha(\mu)$.
Since $\mu - w_\alpha(\mu) = \left<\alpha,\mu \right>\alpha$
is a positive integral multiple of $\alpha$, it follows that
$\lambda + \alpha$ lies on the edge of this triangle with vertices $\lambda$
and $\lambda + \mu - w_\alpha(\mu)$.  Since $K_a$ is a hyperplane
which intersects this side of the triangle, it must also intersect
one of the other sides, i.e., either $K_a$
intersects the line between $\lambda$ and $\lambda+\mu$, or
intersects the line between $\lambda +\mu$ and $\lambda +\mu
- w_\alpha(\mu)$.  The former case satisfies (1).  
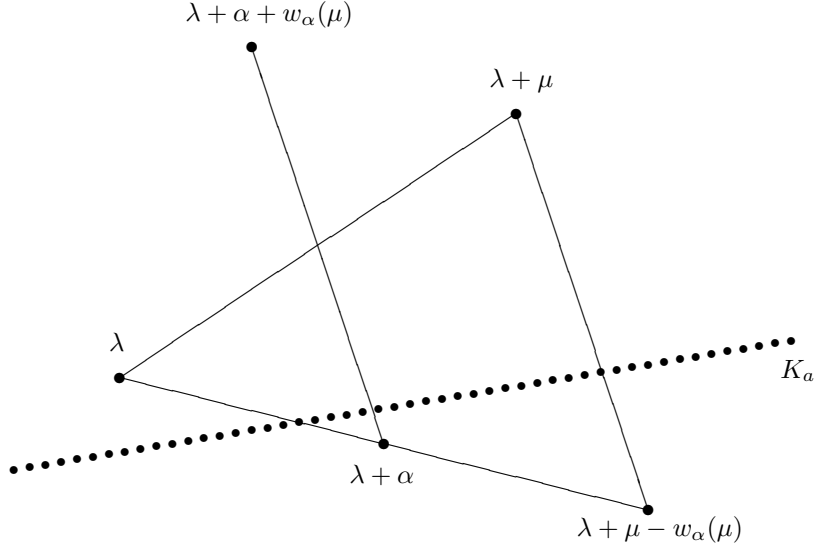
\begin{figure}[htb]
\begin{center}
\begin{picture}(300,250)
\put(50,100){\circle*{4}}
\put(200,200){\circle*{4}}
\put(250,50){\circle*{4}}
\put(50,100){\line(3,2){150}}
\put(50,100){\line(4,-1){200}}
\put(200,200){\line(1,-3){49.5}}
\put(150,75){\circle*{4}}
\put(100,225){\circle*{4}}
\put(150,75){\line(-1,3){49.5}}
\multiput(10,65)(6,1){50}{\circle*{3}}
\put(46,110){$\lambda$}
\put(190,210){$\lambda + \mu$}
\put(223,40){$\lambda + \mu - w_\alpha(\mu)$}
\put(137,60){$\lambda + \alpha$}
\put(75,235){$\lambda + \alpha + w_\alpha(\mu)$}
\put(300,100){$K_a$}
\end{picture}
\end{center}
\caption{The hyperplane $K_a$ and the triangle with vertices $\lambda,
\lambda + \mu$, and $\lambda + \mu - w_\alpha(\mu)$.}
\end{figure}
If the latter
case holds, then as $K_a$ intersects the line between $\lambda$ and
$\lambda + \alpha$ and the line between $\lambda +\mu$ and $\lambda +\mu
- w_\alpha(\mu)$, $K_a$ must also intersect the line between $\lambda +
\alpha$ and $\lambda + \alpha +w_\alpha(\mu)$ (see Figure 1).  
This satisfies (2).
\end{proof}

\begin{lemma}
\label{positivekill}
Suppose $\beta$ is a weight for $\sigma$ such that $\beta \in K_a^+$
and $w_\alpha(\beta) \in K_a^-$.
Then for every $v\in V_\beta\subset \Bbb R^n$ and for every 
maximal weight vector $X$ in our given irreducible subrepresentation 
of $\rho$, $Xv=0$.  Further, this result holds if we replace 
$K_a$ with $K_{w(a)}$ for any $w\in\mathcal{W}$.
\end{lemma}
\begin{proof}
Let $\lambda$ be as in Lemma \ref{killweightlem}, and assume that the
first of the two possible conditions holds.  Since $\lambda\in K_a^-$ and
$\lambda+\mu\in K_a^+$,  Lemma \ref{remB} implies that $Xw=0$ for
any weight vector $w\in \Bbb R^n$ with weight $\lambda$
and for any weight vector $X\in \frak{U}$ with weight $\mu$.
For $v\in V_\beta$, there exists 
a series of root vectors
$\{R_i\}$ with root $\alpha$ and $w\in V_\lambda$ such that 
$(\sigma\pi(R_l))\circ \cdots \circ (\sigma\pi(R_1))(w)=v$.  Since 
$\rho\pi(R_i)(X)=[R_i,X]=0$, we have 
\[ Xv = X(\sigma\pi(R_l))\circ \cdots \circ (\sigma\pi(R_1))(w)
= (\sigma\pi(R_l))\circ \cdots \circ (\sigma\pi(R_1))(Xw) =0. \]

Suppose the latter case of Lemma \ref{killweightlem} holds. 
Choose $g_\alpha \in G$ such that $Ad(g_\alpha)$
corresponds to $w_\alpha$.  
An argument similar to the one above shows that for any
$v\in V_{w_\alpha(\beta)}$ and any maximal $\rho$ weight vector $X$,
$(Ad(g_\alpha)(X))v=0$. 
Note that $\sigma\pi(g_\alpha^{-1})(V_\beta)
\subset V_{w_\alpha(\beta)}$.  Hence, for any $w\in V_\beta$, we can 
write $w=\sigma\pi(g_\alpha)(v)$ for some $v\in V_{w_\alpha(\beta)}$.  
Finally, as $\rho$ comes from the adjoint representation of $H$, we have
$\sigma(\rho(g_\alpha)(X)) = \sigma(g_\alpha)(\sigma(X))\sigma(g_\alpha^{-1})$.
Thus,
\[ 0=(Ad(g_\alpha)(X))v = \sigma(\rho(g_\alpha(X))v 
=\sigma(g_\alpha)\sigma(X)\sigma(g_\alpha^{-1})\sigma(g_\alpha)v
=\sigma(g_\alpha)\sigma(X)v = \sigma(g_\alpha)(Xv). \]
Since $\sigma(g_\alpha)$ is invertible, we have that $Xv=0$.
\end{proof}

The previous lemma shows that certain types of weight spaces are killed
by all weight vectors $X \in\frak{U}$ with weight $\mu$.  
In the remainder of the proof, we first show that
the maximal weight space $V_\chi$ for any irreducible subrepresentation
of $\sigma$ is killed by all such $X$'s, and second, that this
is sufficient to see that all weight spaces for $\sigma$ are killed
by any such $X$.

\begin{prop}
\label{rootprop}
Let $\mathcal{R}$ be a root system with $\mathcal{C}^+$ the 
positive Weyl chamber, $\alpha$ the maximal root, $\mathcal{W}$
the Weyl group, and $w_\alpha$ the element of $\mathcal{W}$ 
defined as reflection through the hyperplane perpendicular to $\alpha$.  
If $\chi \in\overline{\mathcal{C}^+}$
is a weight, then there exists $w \in \mathcal{W}$ such that
$(\chi,w(a))$ and $(w_\alpha(\chi),w(a))$ have opposite signs.
\end{prop}

\begin{proof}
Let $\mathcal{P}$ be the 2 dimensional space spanned by $\chi$ and $\alpha$ 
(see Figure 2).
\begin{figure}[htb]
\begin{center}
\begin{picture}(300,260)(-150,-110)
\put(0,0){\line(2,3){80}}
\put(0,0){\line(-2,-3){80}}
\put(0,0){\line(-5,-3){115}}
\put(0,0){\line(1,6){22}}
\put(0,0){\circle*{5}}
\put(48,72){\circle*{4}}
\put(-48,-72){\circle*{4}}
\put(0,0){\line(-1,2){50}}
\put(0,0){\line(-5,3){90}}
\multiput(0,0)(2,-3){40}{\circle*{3}}
\multiput(0,0)(-2,3){40}{\circle*{3}}
\put(54,69){$\alpha$}
\put(-42,-75){$-\alpha$}
\put(48,-60){$\mathcal{L}$}
\put(46,123){$\mathcal{C}^+$}
\put(-148,-100){$w_\alpha(\mathcal{C}^+)$}
\put(-70,58){$\mathcal{C}_1$}
\put(37,90){\circle*{4}}
\put(-69,-69){\circle*{4}}
\put(-93,-81){$w_\alpha(\chi)$}
\put(34,96){$\chi$}
\end{picture}
\end{center}
\caption{The plane $\mathcal{P}$ spanned by $\chi$ and $\alpha$.}
\end{figure}
Then $\chi \in \mathcal{C}^+ \cap \mathcal{P}$ and 
$w_\alpha(\chi) \in w_\alpha(\mathcal{C}^+)
\cap \mathcal{P}$.  Moreover, since $(\nu,\alpha) > 0$ for 
every $\nu \in\overline{\mathcal{C}^+}$, 
$\overline{\mathcal{C}^+} \cap \mathcal{P}$ and $w_\alpha(\overline{
\mathcal{C}^+}) \cap \mathcal{P}$ cannot share an edge.
Hence, there exists a half cone in $\mathcal{P}$ which 
separates $\mathcal{C}^+ \cap
\mathcal{P}$ and $w_\alpha(\mathcal{C}^+) \cap 
\mathcal{P}$ and which intersects the
line segment joining $\chi$ and $w_\alpha(\chi)$.  Note that such a
half cone consists of the intersections of Weyl chambers with $\mathcal{P}$.
Hence we can find some Weyl chamber $\mathcal{C}_1$ such that $\mathcal{C}_1$
intersects the line segment joining $\chi$ and $w_\alpha(\chi)$.

If $a$ is normally hyperbolic, then so is $\omega(a)$ for every 
$\omega\in\mathcal{W}$.
Thus, using the transitivity of $\mathcal{W}$ on the set of Weyl
chambers, we can find some $w\in\mathcal{W}$ such that $K_{w(a)}$ intersects
$\mathcal{C}_1$ which, since $K_{w(a)}$ has codimension 1, must
intersect $\mathcal{P}$ in at least a line.  
Let $\mathcal{L}= K_{w(a)} \cap \mathcal{P}$.  
If $\mathcal{L}=\mathcal{P}$, then $w(a)$ is perpendicular to $\alpha$.  
This cannot happen since we have chosen $a$ such that
$(\omega(a),\beta) \ne 0$ for every $\omega\in\mathcal{W}$ and 
for every root $\beta$.

We may therefore assume that $\mathcal{L}$ is a line. Since $\mathcal{C}_1$
intersects the line segment between $\chi$ and $w_\alpha(\chi)$, it 
follows that $\mathcal{L}$ does as well.  Hence $(\chi,w(a))$ and
$(w_\alpha(\chi),w(a))$ have opposite signs.  This is possible
since $w(a)$ is conjugate to $a$ in $G$ and hence normally hyperbolic.
\end{proof}

\begin{cor}
If $V_\chi$ is a maximal weight space for any irreducible subrepresentation
of $\sigma$, 
then $Xv=0$ for every $v\in V_\chi$ and for every weight vector
$X \in \frak{U}$ with weight $\mu$.
\end{cor}
\begin{proof}
Proposition \ref{rootprop} shows that, by possibly having to replace
$a$ with $-a$, there 
exists some $w \in \mathcal{W}$ such that 
$\chi \in K^+_{w(a)}$ and 
$w_\alpha(\chi) \in K^-_{w(a)}$.  
Now apply Lemma \ref{positivekill} using $K_{w(a)}^{\pm}$ in place
of $K_a^\pm$.
\end{proof}

Let $V_\nu\subset \frak{U}$ be any weight space in 
$\sigma$ with weight $\nu$.  We will
complete the proof by showing that any $X \in V_\mu$ kills any
$v\in V_\nu$.  Assume for the sake of contradiction that $Xv \ne 0$.
Let $\{\chi_i\}$ be the set of the maximal weights for all the
irreducible subrepresentations for $\sigma$.  From Proposition
\ref{rootprop}, we can choose $a_i\in\mathcal{W}(\pm a)$ 
to be a normally
hyperbolic element in $G$ such that $\chi_i \in K^+_{a_i} \cap
K^-_{w_\alpha(a_i)}$.  If $\nu \in K^+_{a_i} \cap
K^-_{w_\alpha(a_i)}$ for any $i$, then Lemma \ref{positivekill} would 
imply $Xv = 0$.  Hence we can assume that $\nu \not\in \mathcal{K} =
\bigcup_i  (K^+_{a_i} \cap K^-_{w_\alpha(a_i)})$.

Further, if $\nu + \mu \in\mathcal{K}$, then by Lemma \ref{remB}
we would have $Xv=0$. Thus $\nu + \mu \not\in \mathcal{K}$, and 
we have that $Xv$ cannot be a maximal weight vector for any irreducible 
subrepresentation of $\sigma$.  Hence it is possible to find some positive
simple root vector $R_1$ such that $\sigma\pi(R_1)(Xv)\ne 0$.  Let $v_1 =
\sigma\pi(R_1)v$.  Then $v_1$ is a weight vector with some weight $\nu_1$, 
such that with respect to the usual ordering of weights, $\nu \prec \nu_1$.
Since $\rho\pi(R_1)(X) = 0$, we have $Xv_1 = X(\rho\pi(R_1)(v)) =
\sigma\pi(R_1)(Xv) \ne 0$.  We repeat this process using $v_1$ instead of
$v$.  Since $Xv_1 \ne 0$, it follows that 
both $\nu_1$ and $\nu_1 + \mu$ do not lie in $\mathcal{K}$.
Therefore, $Xv_1$ cannot be a maximal weight for 
any irreducible subrepresentation of
$\sigma$.  Hence we can find a positive simple root $R_2$ such that
$\sigma\pi(R_2)(Xv_1)\ne 0$.  If $v_2 = \sigma\pi(R_2)(v_1)$,
then it follows that $Xv_2 \ne 0$, and $v_2$ is a weight vector
with some weight $\nu_2$ such that $\nu\prec\nu_1\prec\nu_2$.  
We can continue this process and 
produce an infinite sequence of weight vectors each with a different weight.
Obviously, this contradicts the finite dimensionality of $\sigma$.
Hence our assumption that $Xv \ne 0$ must be false.  Theorem \ref{THM2}
now follows. \qed

\subsection{Proof of Proposition \ref{THM1}}

As mentioned before the statement of Proposition \ref{THM1}, by setting
$\phi(m) = E^-_a(m)$ we have an $A$-invariant H\"{o}lder section
of $E_{Gr(k,n)} \to M$.  Let $\Phi:P \to Gr(k,n)$ be the $H$-equivariant
map corresponding to $\phi$.

\begin{lemma}
The section $\phi$ is invariant under some parabolic subgroup of $G$.
\end{lemma}
\begin{proof}
If $U^- = \{ u \in G \mid Ad(a^n)u \to 1 \text{ as } n \to \infty\}$, then
$Q = Z_G(A)\cdot U^-$ is a parabolic subgroup.  To see that
$\phi$ is $Q$-invariant note that for a fixed $q \in Q$, $\{a^nqa^{-n}\}$
is bounded in $G$.  Hence, for $v \in E^-_a(m)$, 
$d(a^n)d(q)v = d(a^nqa^{-n})d(a^n)v
\to 0$ as $n \to \infty$ showing that $\phi$ is indeed $Q$-invariant.
\end{proof}

Proposition \ref{THM1} now follows by letting $K$ be the kernel
of the $H$ action on $\Phi(P)$. \qed

\subsection{Proof of Theorem \ref{secthm}}
The Theorem follows immediately by combining Proposition \ref{THM1}
and Theorem \ref{TSR}.  It remains only to see that $K$ is compact
when the irreducible subrepresentations of $\pi$ are multiplicity
free.

Let $H$ be the H\"{o}lder hull, so by Theorem \ref{THM2}, $H$ is reductive
with compact center.  Write $H =CS$ as an almost direct product where
$C$ is compact and $S$ is semisimple with no compact factors. 
Since $G$ has no compact factors, $\pi(G)$ must take values in $S$.
Since $K$ is the kernel of the $H$ action on $\Phi(P)$, it commutes with
the image of $\pi(G)$.  Using Schur's Lemma, and the assumption of 
multiplicity free, $K \cap S$ must be abelian.  Thus, $K \cap S$ has
a trivial connected component, i.e., we have that $K \subset C$. \qed

\bibliography{goetze}

\end{document}